\documentclass[11pt]{article}
\usepackage[dvips]{graphicx}
\usepackage{textcomp}
\usepackage{amsbsy}
\usepackage{latexsym}
\usepackage[mathscr]{eucal}
\usepackage{amsfonts}
\usepackage{amssymb}
\usepackage[usenames]{color}

\usepackage{fullpage}

\begin{document}
\bibliographystyle{plain}
\title
{
A general wavelet-based profile decomposition in the critical embedding of function spaces
\thanks{%
 The third author was supported by the EPSRC Science and Innovation award to the Oxford Centre for Nonlinear PDE (EP/E035027/1).}}
\author{
Hajer Bahouri, Albert Cohen and Gabriel Koch
}
\hbadness=10000
\vbadness=10000
\newtheorem{lemma}{Lemma}[section]
\newtheorem{prop}[lemma]{Proposition}
\newtheorem{cor}[lemma]{Corollary}
\newtheorem{theorem}[lemma]{Theorem}
\newtheorem{remark}[lemma]{Remark}
\newtheorem{example}[lemma]{Example}
\newtheorem{definition}[lemma]{Definition}
\newtheorem{proper}[lemma]{Properties}
\newtheorem{assumption}[lemma]{Assumption}
%
\def\RR{\rm \hbox{I\kern-.2em\hbox{R}}}
\def\NN{\rm \hbox{I\kern-.2em\hbox{N}}}
\def\ZZ{\rm {{\rm Z}\kern-.28em{\rm Z}}}
\def\CC{\rm \hbox{C\kern -.5em {\raise .32ex \hbox{$\scriptscriptstyle
|$}}\kern
-.22em{\raise .6ex \hbox{$\scriptscriptstyle |$}}\kern .4em}}
\def\vp{\varphi}
\def\<{\langle}
\def\>{\rangle}
\def\t{\tilde}
\def\i{\infty}
\def\e{\varepsilon}
\def\sm{\setminus}
\def\nl{\newline}
\def\o{\bar}
\def\wt{\widetilde}
\def\wh{\widehat}
\def\cT{{\cal T}}
\def\cA{{\cal A}}
\def\cI{{\cal I}}
\def\cV{{\cal V}}
\def\cB{{\cal B}}
\def\cR{{\cal R}}
\def\cD{{\cal D}}
\def\cP{{\cal P}}
\def\cJ{{\cal J}}
\def\cM{{\cal M}}
\def\cO{{\cal O}}
\def\Chi{\raise .3ex
\hbox{\large $\chi$}} \def\vp{\varphi}
\def\lsima{\hbox{\kern -.6em\raisebox{-1ex}{$~\stackrel{\textstyle<}{\sim}~$}}\kern -.4em}
\def\lsim{\hbox{\kern -.2em\raisebox{-1ex}{$~\stackrel{\textstyle<}{\sim}~$}}\kern -.2em}
\def\[{\Bigl [}
\def\]{\Bigr ]}
\def\({\Bigl (}
\def\){\Bigr )}
\def\[{\Bigl [}
\def\]{\Bigr ]}
\def\({\Bigl (}
\def\){\Bigr )}
\def\L{\pounds}
\def\pr{{\rm Prob}}
\newcommand{\cs}[1]{{\color{blue}{#1}}}
\def\ds{\displaystyle}
\def\ev#1{\vec{#1}}     
\newcommand{\lt}{\ell^{2}(\nabla)}
\def\Supp#1{{\rm supp\,}{#1}}
\def\R{\mathbb{R}}
\def\E{\mathbb{E}}
\def\nl{\newline}
\def\T{{\relax\ifmmode I\!\!\hspace{-1pt}T\else$I\!\!\hspace{-1pt}T$\fi}}
\def\N{\mathbb{N}}
\def\Z{\mathbb{Z}}
\def\N{\mathbb{N}}
\def\Zd{\Z^d}
\def\Q{\mathbb{Q}}
\def\C{\mathbb{C}}
\def\Rd{\R^d}
\def\gsim{\mathrel{\raisebox{-4pt}{$\stackrel{\textstyle>}{\sim}$}}}
\def\sime{\raisebox{0ex}{$~\stackrel{\textstyle\sim}{=}~$}}
\def\lsim{\raisebox{-1ex}{$~\stackrel{\textstyle<}{\sim}~$}}
\def\div{\mbox{ div }}
\def\M{M}  \def\NN{N}                  
\def\L{{\ell}}               
\def\Le{{\ell^1}}            
\def\Lz{{\ell^2}}
\def\Let{{\tilde\ell^1}}     
\def\Lzt{{\tilde\ell^2}}
\def\Ltw{\ell^\tau^w(\nabla)}
\def\t#1{\tilde{#1}}
\def\la{\lambda}
\def\La{\Lambda}
\def\ga{\gamma}
\def\BV{{\rm BV}}
\def\Ga{\eta}
\def\al{\alpha}
\def\cZ{{\cal Z}}
\def\cA{{\cal A}}
\def\cU{{\cal U}}
\def\cV{{\cal V}}
\def\argmin{\mathop{\rm argmin}}
\def\argmax{\mathop{\rm argmax}}
\def\prob{\mathop{\rm prob}}

\def\cO{{\cal O}}
\def\cA{{\cal A}}
\def\cC{{\cal C}}
\def\cF{{\cal F}}
\def\bu{{\bf u}}
\def\bz{{\bf z}}
\def\bZ{{\bf Z}}
\def\bI{{\bf I}}
\def\cE{{\cal E}}
\def\cD{{\cal D}}
\def\cG{{\cal G}}
\def\cI{{\cal I}}
\def\cJ{{\cal J}}
\def\cM{{\cal M}}
\def\cN{{\cal N}}
\def\cT{{\cal T}}
\def\cU{{\cal U}}
\def\cV{{\cal V}}
\def\cW{{\cal W}}
\def\cL{{\cal L}}
\def\cB{{\cal B}}
\def\cG{{\cal G}}
\def\cK{{\cal K}}
\def\cS{{\cal S}}
\def\cP{{\cal P}}
\def\cQ{{\cal Q}}
\def\cR{{\cal R}}
\def\cU{{\cal U}}
\def\bL{{\bf L}}
\def\bl{{\bf l}}
\def\bK{{\bf K}}
\def\bC{{\bf C}}
\def\X{X\in\{L,R\}}
\def\ph{{\varphi}}
\def\D{{\Delta}}
\def\H{{\cal H}}
\def\bM{{\bf M}}
\def\bx{{\bf x}}
\def\bj{{\bf j}}
\def\bG{{\bf G}}
\def\bP{{\bf P}}
\def\bW{{\bf W}}
\def\bT{{\bf T}}
\def\bV{{\bf V}}
\def\bv{{\bf v}}
\def\bt{{\bf t}}
\def\bz{{\bf z}}
\def\bw{{\bf w}}
\def \span{{\rm span}}
\def \meas {{\rm meas}}
\def\rhom{{\rho^m}}
\def\lll{\langle}
\def\argmin{\mathop{\rm argmin}}
\def\argmax{\mathop{\rm argmax}}
\def\dJ{\nabla}
\newcommand{\ba}{{\bf a}}
\newcommand{\bb}{{\bf b}}
\newcommand{\bc}{{\bf c}}
\newcommand{\bd}{{\bf d}}
\newcommand{\bs}{{\bf s}}
\newcommand{\bff}{{\bf f}}
\newcommand{\bp}{{\bf p}}
\newcommand{\bg}{{\bf g}}
\newcommand{\by}{{\bf y}}
\newcommand{\br}{{\bf r}}
\newcommand{\be}{\begin{equation}}
\newcommand{\ee}{\end{equation}}
\newcommand{\bea}{$$ \begin{array}{lll}}
\newcommand{\eea}{\end{array} $$}
\def \Vol{\mathop{\rm  Vol}}
\def \mes{\mathop{\rm mes}}
\def \Prob{\mathop{\rm  Prob}}
\def \exp{\mathop{\rm    exp}}
\def \sign{\mathop{\rm   sign}}
\def \sp{\mathop{\rm   span}}
\def \vphi{{\varphi}}
\def \csp{\overline \mathop{\rm   span}}
\newcommand{\KL}{Karh\'unen-Lo\`eve }
%
\newcommand{\beqn}{\begin{equation}}
\newcommand{\eeqn}{\end{equation}}
\def\beginproof{\noindent{\bf Proof:}~ }
\def\endproof{\hfill\rule{1.5mm}{1.5mm}\\[2mm]}

\newenvironment{Proof}{\noindent{\bf Proof:}\quad}{\endproof}

\renewcommand{\theequation}{\thesection.\arabic{equation}}
\renewcommand{\thefigure}{\thesection.\arabic{figure}}

\makeatletter
\@addtoreset{equation}{section}
\makeatother

\newcommand\abs[1]{\left|#1\right|}
\newcommand\clos{\mathop{\rm clos}\nolimits}
\newcommand\trunc{\mathop{\rm trunc}\nolimits}
\renewcommand\d{d}
\newcommand\dd{d}
\newcommand\diag{\mathop{\rm diag}}
\newcommand\dist{\mathop{\rm dist}}
\newcommand\diam{\mathop{\rm diam}}
\newcommand\cond{\mathop{\rm cond}\nolimits}
\newcommand\eref[1]{{\rm (\ref{#1})}}
\newcommand{\iref}[1]{{\rm (\ref{#1})}}
\newcommand\Hnorm[1]{\norm{#1}_{H^s([0,1])}}
\def\int{\intop\limits}
\renewcommand\labelenumi{(\roman{enumi})}
\newcommand\lnorm[1]{\norm{#1}_{\ell^2(\Z)}}
\newcommand\Lnorm[1]{\norm{#1}_{L^2([0,1])}}
\newcommand\LR{{L^2(\R)}}
\newcommand\LRnorm[1]{\norm{#1}_\LR}
\newcommand\Matrix[2]{\hphantom{#1}_#2#1}
\newcommand\norm[1]{\left\|#1\right\|}
\newcommand\ogauss[1]{\left\lceil#1\right\rceil}
\newcommand{\QED}{\hfill
\raisebox{-2pt}{\rule{5.6pt}{8pt}\rule{4pt}{0pt}}%
  \smallskip\par}
\newcommand\Rscalar[1]{\scalar{#1}_\R}
\newcommand\scalar[1]{\left(#1\right)}
\newcommand\Scalar[1]{\scalar{#1}_{[0,1]}}
\newcommand\Span{\mathop{\rm span}}
\newcommand\supp{\mathop{\rm supp}}
\newcommand\ugauss[1]{\left\lfloor#1\right\rfloor}
\newcommand\with{\, : \,}
\newcommand\Null{{\bf 0}}
\newcommand\bA{{\bf A}}
\newcommand\bB{{\bf B}}
\newcommand\bR{{\bf R}}
\newcommand\bD{{\bf D}}
\newcommand\bE{{\bf E}}
\newcommand\bF{{\bf F}}
\newcommand\bH{{\bf H}}
\newcommand\bU{{\bf U}}
\newcommand\cH{{\cal H}}
\newcommand\sinc{{\rm sinc}}
\def\enorm#1{| \! | \! | #1 | \! | \! |}

\newcommand{\dm}{\frac{d-1}{d}}

\let\bm\bf
\newcommand{\bbeta}{{\mbox{\boldmath$\beta$}}}
\newcommand{\bal}{{\mbox{\boldmath$\alpha$}}}
\newcommand{\bbi}{{\bm i}}

\def\nnew{\color{Red}}
\def\mnew{\color{Blue}}

\newcommand{\dI}{\Delta}
\maketitle
\date{}

\begin{abstract}
We characterize the lack of compactness in the critical embedding
of functions spaces $X\subset Y$ having similar scaling properties
in the following terms : a sequence $(u_n)_{n\geq 0}$ bounded in $X$
has a subsequence that can be expressed as a finite sum of translations
and dilations of functions $(\phi_l)_{l>0}$ such that the remainder converges to zero
in $Y$ as the number of functions in the sum and $n$ tend to $+\infty$.
Such a decomposition was established by G\'erard in \cite{Ge}
for the embedding of the homogeneous Sobolev space $X=\dot H^s$
into the $Y=L^p$ in $d$ dimensions with $0<s=d/2-d/p$,
and then generalized by Jaffard in \cite{Ja} to the case
where $X$ is a Riesz potential space,
using wavelet expansions. In this paper, we revisit the
wavelet-based profile decomposition, in order to
treat a larger range of examples of critical embedding in a hopefully
simplified way. In particular we identify
two generic properties on the spaces $X$ and $Y$
that are of key use in building the profile decomposition.
These properties may then easily be checked for typical
choices of $X$ and $Y$ satisfying critical embedding properties.
These includes Sobolev, Besov, Triebel-Lizorkin, Lorentz, H\"older and BMO spaces.
\end{abstract}

\section{Introduction}

The critical embedding of homogeneous Sobolev spaces
in dimension $d$ states that for $0\leq t<s$
and $1\leq p<q<\infty$ such that $d/p-d/q=s-t$, one has
\be
\dot W^{s,p}(\R^d)\subset \dot W^{t,q}(\R^d).
\ee
The lack of compactness in this embedding
can be described
in terms of an asymptotic decomposition
following G\'erard \cite{Ge} who considered the case
$p=2$ and $t=0$, and Jaffard \cite{Ja} who considered general values
$p>1$ with the Riesz potential spaces $\dot H^{s,p}$ in replacement
of Sobolev spaces $\dot W^{s,p}$, again with $t=0$.
Their results can be formulated in the following terms :
a sequence $(u_n)_{n\geq 0}$ bounded in $\dot H^{s,p}(\RR^d)$
can be decomposed up to a subsequence extraction according to
\be
u_n=\sum_{l=1}^L h_{l,n}^{s-d/p}\phi^l\(\frac {\cdot - x_{l,n}}{h_{l,n}}\) + r_{n,L}
\label{decprof}
\ee
where $(\phi^l)_{l>0}$ is a family of functions
in $\dot H^{s,p}(\RR^d)$ and
where
$$
\lim_{L\to +\infty}\(\limsup_{n\to +\infty}\|r_{n,L}\|_{L^{q}}\)=0.
$$
This
decomposition is ``asymptotically orthogonal'' in the sense that for $k\neq l$
$$
|\log (h_{l,n}/h_{k,n})| \to +\infty
\;\;{\rm or}
\;\;
|x_{l,n}-x_{k,n}|/h_{l,n}\to +\infty,
\;\;{\rm as}\;\; n\to +\infty.
$$
This type of decomposition was also obtained earlier
in \cite{BC} for a bounded sequence in $H_0^{1}(D, \R^3)$ of solutions of an elliptic problem, with $D$
the open unit disk of $\R^2$ and in \cite{St} and \cite{So} for the critical injections
of $W^{1,2}(\Omega)$ in Lebesgue space
and of $W^{1,p}(\Omega)$ in Lorentz spaces
respectively, with $\Omega$ a bounded domain of $\R^d$.
They were also studied in \cite{ST} in an abstract Hilbert space
framework and in \cite{Benameur} in  the Heisenberg group context.\\

The above mentioned references treat different types of examples of critical
embedding by different methods. One of the motivations of
the present paper is to identify some fundamental mechanisms that
lead to such results for a general critical embedding
$$
X \subset Y,
$$
in a unified way. Here $X$ and $Y$ are generic homogeneous function spaces which, similar to
the above particular cases, have the same scaling properties
in the sense that for any function $f$ and $h>0$
\be
\|f(h\cdot)\|_X= h^r \|f\|_X \;\;{\rm and}\;\; \|f(h\cdot)\|_Y= h^r \|f\|_Y,
\label{scaling}
\ee
for the same value of $r$.
In a similar way to Jaffard, we use wavelet bases in
order to construct the functions $\phi_l$,
yet in a somehow different and hopefully simpler way.
Our construction is based on two basic key properties
of wavelet expansions in the spaces $X$ and $Y$,
which may then be easily checked on particular pairs
of spaces of interest. In particular, any critical embedding involving Sobolev, Besov, Triebel-Lizorkin,
Lorentz, H\"older or BMO spaces is covered by our approach.\\

The study of the lack of compactness in the
critical embedding of Sobolev spaces supplies us with a large amount of information about
solutions of nonlinear  partial differential equations, both in the
elliptic frame or the evolution frame. One has, for example,
\begin{itemize}
\item
the pioneering works of P. -L.  Lions \cite{Lions1} and
\cite{Lions-II} for the sake of
geometric problems,
\item
the
description of bounded energy sequences of solutions to the
defocusing semi-linear quintic wave equation, up to remainder terms small in energy norm in
\cite{BG},
\item the characterization
 of the defect of compactness for  Strichartz estimates  for the Shr\"odinger equation in \cite{ker},
\item
the understanding of features of solutions
of nonlinear wave equations with exponential growth in \cite{BMM},
\item
the sharp estimate of the time life-span of the
focusing critical semi-linear  wave equation by means of the size of
energy of the Cauchy data in \cite{km},
\item
the study of the bilinear Strichartz estimates for the wave equation in \cite{Tao}.
\end{itemize}
For further applications, we refer to  \cite{GG}, \cite{Ga}, \cite{IB}, \cite{Ma}, \cite{La} and the references therein.\\

Our results which cover  a broad spectrum of spaces could be at the origin
of several prospectus of similar types of regularity results for Navier-Stokes systems (as in \cite{kk,gkp}),
qualitative study of non linear evolution equations or estimates of the span life of focusing semi-linear dispersive evolution equations.

\subsection{Wavelet expansions}

Wavelet decompositions of a function have the form
\be
f=\sum_{\lambda\in \nabla} d_\lambda\psi_\lambda,
\label{decomp}
\ee
where $\lambda=(j,k)$ concatenates the scale index $j=j(\lambda)$ and space index
$k=k(\lambda)$: for $d=1$, we have with the $L^2$ normalization,
$$
\psi_{j,k}=\psi_\lambda=2^{j/2}\psi (2^j\cdot-k), \;\; j\in\Z, \; k\in \Z,
$$
where $\psi$ is the so-called ``mother wavelet''. In higher dimension $d>1$,
one needs several generating functions $\psi^e$ for $e\in E$ a finite set, so that
setting $\psi_\lambda:=(\psi^e_\lambda)^T_{e\in E}$ and $d_\lambda=(d_\lambda^e)_{e\in E}$,
we can again write \iref{decomp} with $d_\lambda\psi_\lambda$ a finite
dimensional inner product and
$$
\psi_{j,k}=\psi_\lambda=2^{dj/2}\psi (2^j\cdot-k), \;\; j\in\Z, \; k\in\Z^d.
$$
The index set $\nabla$ in \iref{decomp} is thus always defined as
$$
\nabla:=\Z\times \Z^d
$$
Note that $\lambda$ may also be identified to a dyadic cube
$$
\lambda\sim 2^{-j}(k+[0,1]^d).
$$
We shall sometimes use the notation
$$
|\lambda|:=j(\lambda),
$$
for the scale level of $\lambda$.
In all the sequel, we systematically normalize our wavelets in
$X$ which is equivalent to normalizing them in $Y$ in view of \iref{scaling}:
\be
\psi_{j,k}=\psi_\lambda=2^{rj}\psi (2^j\cdot-k).
\label{norm}
\ee
It is known that, in addition to being Schauder bases, wavelet bases are unconditional bases
for ``most'' classical function spaces, including
in particular the family of Besov and Triebel-Lizorkin spaces:
for such spaces $X$ there exists a constant $D$ such that
for any finite subset $E\subset \nabla$ and coefficients
vectors $(c_\lambda)_{\lambda\in E}$ and $(d_\lambda)_{\lambda\in E}$
such that $|c_\lambda| \leq |d_\lambda|$ for all $\lambda$, one has
\be
\|\sum_{\lambda\in E} c_\lambda\psi_\lambda\|_X \leq D\|\sum_{\lambda\in E} d_\lambda\psi_\lambda\|_X.
\label{uncond}
\ee
We refer to \cite{Co,Dau,Me} for more details on the construction
of wavelet bases and on the characterization of classical
function spaces by expansions in such bases.

\subsection{Main results}

Our profile decomposition relies on two key assumptions concerning
wavelet decompositions and the spaces $X$ and $Y$.
\nl

{\it In addition we always work under the general assumption that
our wavelet basis $(\psi_\lambda)_{\lambda\in\Lambda}$ is an unconditional basis
for both spaces $X$ and $Y$. We therefore
assume that \iref{uncond} holds
with some constant $D$ for both norms.}
\nl

Our first assumption
involves the {\it nonlinear projector} that we define for each
$M>0$ as follows: if $f\in X$ has the expansion in the wavelet basis
given by \iref{decomp}, then
\be
Q_M f:=\sum_{\lambda \in E_M} d_\lambda\psi_\lambda,
\label{nonlinear}
\ee
where $E_M=E_M(f)$
is the subset of $\nabla$ of cardinality $M$ that corresponds to the $M$ largest
values of $|d_\lambda|$.

Such a set always exists due to the fact that  $(\psi_\lambda)_{\lambda\in\Lambda}$
is a Schauder basis for $X$, since this implies that for any $\eta>0$
only finitely many coefficients $d_\lambda$ are larger than $\eta$ in modulus.
This set may however not be unique when
some $|d_\lambda|$ are equal, in which case we
may choose an arbitrary realization of such a set.
Recall that we have assumed the normalization
\iref{norm} making $\|\psi_\lambda\|_X$ or $\|\psi_\lambda\|_Y$ independent of $\lambda$,
therefore $E_M$ also corresponds to the $M$ largest
$\|d_\lambda\psi_\lambda\|_X$ or $\|d_\lambda\psi_\lambda\|_Y$.
\nl
\nl
{\bf Assumption 1:} {\it The nonlinear projection satisfies}
\be
\lim_{M\to +\infty}\max_{\|f\|_X \leq 1}\|f-Q_M f\|_Y=0.
\label{assumption1}
\ee
\nl
The fact that the convergence of $Q_Mf$ towards $f$ in
$Y$ holds uniformly on the unit ball of $X$
is tied to the nonlinear nature of the operator $Q_M$: if instead
we took $Q_M$ to be the projection onto a {\it fixed} $M$-dimensional space,
then \iref{assumption1} would be in contradiction with the fact that the
critical embedding of $X$ into $Y$ is not compact. As will be
recalled further, nonlinear approximation
theory actually allows for a more precise quantification of the above property
in most cases of interest, through an estimate of the form
$$
\max_{\|f\|_X \leq1}\|f-Q_M f\|_Y\leq CM^{-s}, \;\; M>0
$$
for some $s>0$ and $C$ only depending on the choice of $X$ and $Y$. However,  Assumption 1 alone will
be sufficient for our purpose.
\nl

Our second assumption only concerns the behavior of
wavelet expansions with respect to the $X$ norm. It reflects
the fact that this norm is stable with respect to certain operations
such as ``shifting'' the indices of wavelet coefficients, as well
as perturbating the value of these coefficients. This is expressed
as follows.
\nl
\nl
{\bf Assumption 2:} {\it Consider a sequence of functions $(f_n)_{n>0}$ which are uniformly bounded in
$X$ and may be written as
\be
f_n=\sum_{\lambda\in\nabla} c_{\lambda,n}\psi_{\lambda},
\ee
and such that for all $\lambda$, the sequence $c_{\lambda,n}$ converges towards a finite limit $c_\lambda$ as $n\to +\infty$.
Then, the series $\sum_{\lambda\in\nabla} c_\lambda\psi_{\lambda}$ converges in $X$ with
\be
\|\sum_{\lambda\in\nabla} c_\lambda\psi_{\lambda}\|_X \leq C \liminf_{n\to +\infty}\|f_n\|_X,
\label{fatou}
\ee
where $C$ is a constant only depending on the space $X$ and on the choice of the wavelet basis.}
\nl
\nl
As will be recalled further, for practical choices of
$X$ such as Besov or Triebel-Lizorkin spaces, the $X$ norm
of a function is equivalent to the norm
of its wavelet coefficients in a certain sequence space. This
allows us to establish \iref{fatou} essentially by invoking Fatou's lemma.

We are now in position to state the main theorem of this paper.
For any function $\phi$, not necessarily a wavelet, and any scale-space index
$\lambda=(j,k)$ we use the notation
\be
\phi_\lambda:=2^{rj}\phi(2^j\cdot-k),
\label{rescale}
\ee
for the version of $\phi$ scaled and translated according to $\lambda$.

\begin{theorem}
\label{maintheo}
Assume that $X$ and $Y$ are two function spaces with
the same scaling \iref{scaling} and continuous embedding $X\subset Y$,
and assume that there exists a wavelet basis $(\psi_\lambda)_{\lambda\in\nabla}$ which is unconditional for both $X$ and $Y$,
and such that Assumptions 1 and 2 hold.
Let $(u_n)_{n>0}$ be a bounded sequence in $X$. Then, up to
subsequence extraction, there exists a family of functions $(\phi^l)_{l>0}$ in $X$
and sequences of scale-space indices $(\lambda_l(n))_{n>0}$ for each $l>0$
such that
\be
u_n=\sum_{l=1}^L \phi_{\lambda_l(n)}^l+ r_{n,L},
\label{profiledec}
\ee
where
$$
\lim_{L\to +\infty}\(\limsup_{n\to +\infty}\|r_{n,L}\|_{Y}\)=0.
$$
The decomposition \iref{profiledec}
is asymptotically orthogonal in the sense that  for any $k\neq l$,
\be
|j(\lambda_k(n))-j(\lambda_l(n))| \to +\infty\;\;{\rm or}\;\; |k(\lambda_k(n))-2^{j(\lambda_k(n))-j(\lambda_l(n))}k(\lambda_l(n))| \to +\infty,
\;\;{\rm as}\;\; n\to +\infty.
\label{orthoscale}
\ee
\end{theorem}

Moreover, we have the following for the specific case where $X$ is a Besov or Triebel-Lizorkin space:
\begin{theorem}\label{genstab}
The decomposition in Theorem \ref{maintheo} is stable in the sense that, for some $\tau = \tau(X)$ we have
\be
\| (\|\phi^l\|_X)_{l>0} \|_{\ell^\tau}\leq CK.
\label{stabTriebel}
\ee
where $C$ is a constant that only depends on $X$ and on the
choice of the wavelet basis and where $K:=\sup_{n\geq 0} \|u_n\|_X$.
\end{theorem}

\begin{remark}
For certain sequences $(u_n)_{n>0}$, it is possible that for any $L>0$ the
decomposition \iref{profiledec} only involves a finite number of profiles
$\phi^l$ for $l=1,\cdots,L_0$, which means that $\phi^l=0$ for
$l>L_0$. Inspection of our proof shows that
the theorem remains valid in such a case,
in the sense that
\be
u_n=\sum_{l=1}^{L_0} \phi_{\lambda_l(n)}^l+ r_{n},
\label{profiledecfinite}
\ee
where
$$
\lim_{n\to +\infty}\|r_{n}\|_{Y}=0.
$$
In particular, the sequence $(u_n)_{n>0}$ is compact in $Y$ if and only if
$\phi^l=0$ for all $l>0$.
\end{remark}

\begin{remark}
\label{remarknonlin}
Inspection of our proof also shows that
in Assumption 1, we may use for $Q_M$
a more general nonlinear projector than
the one obtained by taking the $M$ largest
values of $|d_\lambda|$. Generally speaking,
we may consider a nonlinear projector
$Q_M$ that has the general form \iref{nonlinear},
where the sets $E_M=E_M(f)$ of cardinality $M$ depend on $f$ and
satisfy
$$
E_M(f)\subset E_{M+1}(f).
$$
Such a generalization appears to be useful when treating
certain types of embedding, see \S 3.
\end{remark}

\subsection{Layout}

The effective construction of the decomposition
is addressed in Section \S 2, together with the proof of Theorem \ref{maintheo}.

In Section \S 3, we discuss examples
of $X$ and $Y$ with critical embedding for which Assumptions 1 and 2
can be proved. This includes all previously treated cases, and many
others such as the embedding of Sobolev, Besov and Triebel-Lizorkin spaces
into Lebesgue, Lorentz, BMO and H\"older spaces, or into other
Sobolev, Besov and Triebel-Lizorkin spaces.

Finally, in \S 4, we prove the stability
Theorem \ref{genstab} for both setting of
Besov and Triebel-Lizorkin spaces.

\section{Construction of the decomposition and proof of Theorem \ref{maintheo}}

In this section, we place ourselves under the assumptions of Theorem \ref{maintheo}.
Let $(u_n)_{n>0}$ be a bounded sequence in the space $X$ and define
$$
K:=\sup_{n>0}\|u_n\|_X<+\infty.
$$
The decomposition construction  and the proof of Theorem \ref{maintheo} proceed
in several steps.
\nl
\nl
{\bf Step 1: rearrangements.} We first introduce the
wavelet decompositions
\be
u_n=\sum_{\lambda\in\nabla} d_{\lambda,n} \psi_\lambda.
\ee
For each $n>0$, we consider the non-increasing rearrangement
$(d_{m,n})_{m>0}$ of $(d_{\lambda,n})_{\lambda\in\nabla}$
according to their moduli. We may therefore
write
\be
u_n=\sum_{m>0} d_{m,n} \psi_{\lambda(m,n)}.
\ee
Using the nonlinear projector $Q_M$ defined by \iref{nonlinear}, we further split this expansion into
\be
u_n=\sum_{m=1}^M d_{m,n} \psi_{\lambda(m,n)} + R_M u_n,
\ee
with $R_M u_n =  u_n  - Q_M u_n$. Combining Assumption 1 with the boundedness of $(u_n)_{n>0}$ in $X$, we infer that
\be
\label{nonlinconv}
\lim_{M\to +\infty} \sup_{n>0}\|R_M u_n\|_Y=0.
\ee
Our next observation is that
if $(\psi_\lambda)_{\lambda\in\nabla}$ is an unconditional basis of $X$
then the coefficients $d_{m,n}$ are uniformly bounded: indeed,
\iref{uncond} implies that the rank one projectors
$$
P_\mu: f=\sum_{\lambda\in\Lambda} d_\lambda\psi_\lambda \mapsto  P_\mu f:=d_\mu\psi_\mu,
$$
satisfy the uniform bound
$$
\|P_\mu\|_{X\to X}\leq D,\;\; \mu\in\nabla.
$$
Since we have assumed that our wavelets are normalized in $X$, for example according
to $\|\psi_\mu\|_X=1$ for all $\mu\in\nabla$, we thus have
$$
\sup_{\lambda,n}|d_{\lambda,n}|=\sup_{m,n}|d_{m,n}| \leq DK.
$$
Up to a diagonal subsequence extraction procedure in $n$, we may therefore assume that
for all $m>0$, the sequence $(d_{m,n})_{n>0}$ converges towards a finite limit that depends on $m$,
$$
d_m=\lim_{n\to +\infty} d_{m,n}.
$$
Note that $(|d_m|)_{m>0}$ is a non-increasing sequence
since all sequences $(|d_{m,n}|)_{m>0}$ are non-increasing. We may thus write
$$
u_n=\sum_{m=1}^M d_{m} \psi_{\lambda(m,n)}+t_{n,M},
$$
where
$$
t_{n,M}:=\sum_{m=1}^M(d_{m,n}- d_m)\psi_{\lambda(m,n)}+ R_M u_n.
$$
\nl
{\bf Step 2: construction of approximate profiles.}
We construct the profiles $\phi^l$ as limit of
sequences $\phi^{l,i}$ obtained by the following algorithm. At the first iteration $i=1$, we set
\be
\phi^{1,1}=d_1\psi,\;\; \lambda_1(n):=\lambda(1,n),\;\; \vp_1(n):=n.
\ee
Assume that after iteration $i-1$, we have constructed $L-1$ functions
$(\phi^{1,i},\cdots,\phi^{L-1,i})$ and scale-space index sequences $(\lambda_1(n),\cdots,\lambda_{L-1}(n))$
with $L\leq i$,
as well as an increasing sequence of positive integers $\vp_{i-1}(n)$ such that
$$
\sum_{m=1}^{i-1}d_m\psi_{\lambda(m,\vp_{i-1}(n))}
= \sum_{l=1}^{L-1}\phi^{l,i}_{\lambda_l(\vp_{i-1}(n))}.
$$
At iteration $i$ we shall use the $i$-th
component $d_i\psi_{\lambda(i,\vp_{i-1}(n))}$ to
either modify one of these functions or build a new one according to the
following dichotomy.

(i) First case: assume that we can extract $\vp_i(n)$ from $\vp_{i-1}(n)$
such that for $l=1,\cdots,L-1$ at least one of the following holds:
\be
\lim_{n\to +\infty} |j(\lambda_l(\vp_{i}(n)))-j(\lambda(i,\vp_{i}(n)))|=+\infty,
\label{scaleortho}
\ee
or
\be
\lim_{n\to +\infty}|k(\lambda(i,\vp_{i}(n)))-2^{j(\lambda(i,\vp_{i}(n)))-j(\lambda_l(\vp_{i}(n)))}
k(\lambda_l(\vp_{i}(n)))| =+\infty.
\label{spaceortho}
\ee
In such a case, we create a new profile and scale-space index sequence by defining
$$
\phi^{L,i}:=d_i\psi,\;\; \lambda_L(n):=\lambda(i,n),
$$
and we set $\phi^{l,i}=\phi^{l,i-1}$ for $l=1,\cdots,L-1$.

(ii) Second case: assume that for some subsequence $\vp_{i}(n)$
of $\vp_{i-1}(n)$ and some $l\in \{1,\cdots,L-1\}$
both \iref{scaleortho}
and \iref{spaceortho} do not hold. Then
it is easy checked that  $j(\lambda_l(\vp_{i}(n)))-j(\lambda(i,\vp_{i}(n)))$
and $k(\lambda(i,\vp_{i}(n)))-2^{j(\lambda(i,\vp_{i}(n)))-j(\lambda_l(\vp_{i}(n)))}
k(\lambda_l(\vp_{i}(n)))$ only take a finite number of values as $n$ varies.
Therefore, up to an additional subsequence extraction, we
may assume that there exists numbers $a$ and $b$ such that for all $n>0$,
\be
j(\lambda(i,\vp_{i}(n)))- j(\lambda_l(\vp_{i}(n)))=a,
\label{scaleclose}
\ee
and
\be
k(\lambda(i,\vp_{i}(n)))-2^{j(\lambda(i,\vp_{i}(n)))-j(\lambda_l(\vp_{i}(n)))}
k(\lambda_l(\vp_{i}(n)))=b.
\label{spaceclose}
\ee
We then update the function $\phi^{l,i-1}$ according to
\be
\phi^{l,i}=\phi^{l,i-1}+d_i 2^{ar}\psi(2^{a}\cdot-b).
\ee
and $\phi^{l',i}=\phi^{l',i-1}$ for $l'\in \{1,\cdots,L-1\}$ and $l'\neq l$.

From this construction, and after extracting a diagonal subsequence which eventually coincides with a subsequence of $\varphi_i(n)$ for each $i$, we see that for each value of $M$ there exists
$L=L(M)\leq M$ such that
$$
\sum_{m=1}^M d_{m} \psi_{\lambda(m,n)}=\sum_{l=1}^L \phi^{l,M}_{\lambda_l(n)}.
$$
More precisely, for each $l=1,\cdots,L$, we have
$$
\phi^{l,M}_{\lambda_l(n)}=\sum_{m\in E(l,M)}d_{m} \psi_{\lambda(m,n)},
$$
where the sets $E(l,M)$ for $i=1,\cdots,L$ constitute
a disjoint partition of $\{1,\cdots,M\}$. Note that $E(l,M)\subset E(l,M+1)$
with $\#(E(l,M+1))\leq \#(E(l,M))+1$. Similarly, the number of profiles
$L(M)$ grows at most by $1$ as we move from $M$ to $M+1$.
As explained in Remark 1.2, it is possible that
$L(M)$ terminates at some maximal value $L_0$. Finally, note
that for any $m,m'\in E_{l,M}$ we have that
\be
j(\lambda(m,n))-j(\lambda(m',n))=a(m,m'),
\label{scaleclose2}
\ee
and
\be
k(\lambda(m,n))-2^{j(\lambda(m,n))-j(\lambda(m',n))}k(\lambda(m',n))=b(m,m'),
\label{spaceclose2}
\ee
where $a(m,m')$ and $b(m,m')$ do not depend on $n$.
\nl
\nl
{\bf Step 3: construction of the exact profiles.} We now
want to define the functions $\phi^{l}$ as the limits in $X$
of $\phi^{l,M}$ as $M\to +\infty$. For this purpose,
we shall make use of Assumption 2, combined with
the scaling property \iref{scaling} of the $X$ norm
and the fact that $(\psi_\lambda)_{\lambda\in \nabla}$
is an unconditional basis. For some
fixed $l$ and $M$ such that $l\leq L(M)$, let us define the functions
$$
g^{l,M}:=\sum_{m\in E(l,M)}d_{m} \psi_{\lambda(m)},
$$
$$
f^{l,M,n}:=\sum_{m\in E(l,M)}d_{m,n} \psi_{\lambda(m)},
$$
with $\lambda(m):=\lambda(m,1)$. From the scaling property \iref{scaling}
and the properties \iref{scaleclose2} and \iref{spaceclose2}, we find that
$$
\|f^{l,M,n}\|_X=\|\sum_{m\in E(l,M)}d_{m,n} \psi_{\lambda(m,n)}\|_X.
$$
Since $\sum_{m\in E(l,M)}d_{m,n} \psi_{\lambda(m,n)}$ is a part of the
expansion of $u_n$, we thus find
that
$$
\|f^{l,M,n}\|_X\leq DK,
$$
where $D$ is the constant in \iref{uncond} and $K:=\sup_{n>0} \|u_n\|_X$.
Invoking Assumption 2, we therefore find that
$g^{l,M}$ converges in $X$ towards a limit $g^l$ as
$M\to +\infty$. We finally notice that, by construction, the $g^{l,M}$ are
rescaled versions of the $\phi^{l,M}$: there exists $A>0$ and $B\in\RR^d$ such that
$$
\phi^{l,M}=2^{Ar}g^{l,M}(2^A\cdot  - B).
$$
By \iref{scaling}, we therefore conclude that $\phi^{l,M}$ converges in $X$ towards a limit $\phi^l:=2^{Ar}g^{l}(2^A\cdot - B)$ as
$M\to +\infty$.
\nl
\nl
{\bf Step 4: conclusion of the proof.} For any given $L>0$, we
may write
$$
u_n=\sum_{l=1}^{L} \phi^{l}_{\lambda_l(n)}+r_{n,L},
$$
where, for any value of $M$ such that $L\leq L(M)$, the remainder $r_{n,L}$ may be decomposed into
\be
\sum_{l=1}^{L}( \phi^{l,M}_{\lambda_l(n)}-\phi^{l}_{\lambda_l(n)})
+\sum_{l=1}^{L}\sum_{m\in E(l,M)}(d_{m,n}-d_m)\psi_{\lambda(m,n)}
+ \sum_{l=L+1}^{L(M)}\sum_{m\in E(l,M)}d_{m,n}\psi_{\lambda(m,n)} +R_M u_n.
\label{rnl}
\ee
Note that each of these terms depend on the chosen value
of $M$ but their sum $r_{n,L}$ is actually independent of $M$.
We rewrite this decomposition as
$$
r_{n,L}=r_1(n,L,M)+r_2(n,L,M),
$$
where $r_1$ and $r_2$ stand for the first and last two terms in \iref{rnl}, respectively. By construction,
all values of $m$ which appear in the third term of \iref{rnl} are between $L+1$ and $M$. Therefore
the last two terms in \iref{rnl} may be viewed as a partial sum of
$$
R_L u_n=\sum_{m>L}d_{m,n}\psi_{\lambda(m,n)}.
$$
Since we have assumed that $(\psi_\lambda)_{\lambda\in\nabla}$ is an unconditional basis for $Y$,
we may therefore write
$$
\| r_2(n,L,M) \|_Y \leq D\|R_L u_n\|_Y.
$$
According to Assumption 1, which is expressed by \iref{nonlinconv},
the right hand side converges to $0$ as $L\to +\infty$ uniformly in $n$ and therefore
$$
\lim_{L\to +\infty} \sup_{n,M}\|r_2(n,L,M) \|_Y=0.
$$
We now consider the first two terms inÊ\iref{rnl}. For the first term, we have
$$
\left \|\sum_{l=1}^{L}( \phi^{l,M}_{\lambda_l(n)}-\phi^{l}_{\lambda_l(n)}) \right \|_X \leq
\sum_{l=1}^{L}\| \phi^{l,M}_{\lambda_l(n)}-\phi^{l}_{\lambda_l(n)}\|_X
=\sum_{l=1}^{L}\| \phi^{l,M}-\phi^{l}\|_X .
$$
Therefore, for any fixed $L$, this term goes to $0$ in $X$ as $M\to +\infty$.
For the second term, we first notice that for any fixed $L$ and $M$
such that  $L\leq L(M)$, all values of $m$ which appear in this term are
less or equal to $M$. Since
we have assumed that $(\psi_\lambda)_{\lambda\in\nabla}$ is an unconditional basis for $X$,
it follows that
$$
\left \|\sum_{l=1}^{L}\sum_{m\in E(l,M)}(d_{m,n}-d_m)\psi_{\lambda(m,n)} \right \|_X \leq
D\left \|\sum_{m=1}^{M}(d_{m,n}-d_m)\psi_{\lambda(m,n)} \right \|_X
\leq CD\sum_{m=1}^{M} |d_{m,n}-d_m|,
$$
where $C=\|\psi\|_X=\|\psi_\lambda\|_X$ for all $\lambda\in\nabla$. Therefore
for any any fixed $L$ and $M$
such that  $L\leq L(M)$, this term goes to $0$ in $X$ as $n\to +\infty$. Combining these observations,
we find that for any fixed $L$ and any $\e>0$, there exists $M$ and $n_0$ such that for
all $n\geq n_0$,
$$
\|r_1(n,L,M) \|_X \leq \e.
$$
By continuous embedding, the same holds for $\|r_1(n,L,M) \|_Y$.
Since $M$ was arbitrary in the decomposition \iref{rnl} of $r_{n,L}$, we obtain that
$$
\lim_{L\to +\infty} \(\limsup_{n\to +\infty}\|r_{n,L}\|_Y\)=0,
$$
which concludes the proof of the theorem.

\section{Examples}

Our main result applies to a large range of critical embedding.
Specifically, we consider
\begin{enumerate}
\item
For the space $X$: spaces of Besov type $\dot B^{s}_{p,a}$ or Triebel-Lizorkin type $\dot F^s_{p,a}$ with $1\leq p<\infty$
and $1\leq a \leq \infty$.
\item
For the space $Y$: spaces of Besov type $\dot B^{t}_{q,b}$, Triebel-Lizorkin type $\dot F^t_{q,b}$, Lebesgue type
$L^q$,
Lorentz type $L^{q,b}$, and the space $BMO$, with $1\leq q\leq \infty$ and $1\leq b \leq \infty$.
\end{enumerate}
Note that Lebesgue spaces may be thought of as a particular case of Triebel-Lizorkin spaces
since $L^q=\dot F^0_{q,2}$, yet we treat them separetely
since several results that we invoke further have been
proved in an isolated manner for the specific case of Lebesgue spaces.

The critical embedding for such spaces imposes that $t<s$ together with the scaling
\be
r=\frac d p-s=\frac d q-t,
\label{rscale}
\ee
where $t=0$ if $Y$ is of Lebesgue or Lorentz type, and $t=0$ and $q=\infty$ if $Y=BMO$.
It also imposes some relations between the fine tuning indices $a$ and $b$. For example for $s>0$ and
$p,q$ such that $\frac d p-s=\frac d q$, the space
$\dot B^{s}_{p,a}$ embeds continuously into $L^{q,b}$ if $b\geq a$.

Note that for non integer $t>0$ the H\"older space $\dot C^t$ coincides with the Besov space
$\dot B^{t}_{\infty,\infty}$, and that for all integer $m\geq 0$, the Sobolev space $\dot W^{m,p}$ coincides with
the Triebel-Lizorkin space $\dot F^m_{p,2}$ when $1<p<\infty$. In particular $L^p= \dot F^0_{p,2}$ for $1<p<\infty$.
For $p=1$, it is known that $\dot F^0_{1,2}$ coincides with the Hardy space ${\cal H}^1$ which is a closed subspace of $L^1$.
We refer to \cite{adams} and \cite{triebel} for an introduction to all such spaces.

It is known that properly constructed wavelet bases
are unconditional for all such spaces, see in particular \cite{Me}. In addition, Besov and Triebel-Lizorkin spaces,
as well as $BMO$, may be
characterized by simple properties on wavelet coefficients. More precisely,
for $f=\sum_{\lambda\in\nabla} d_\lambda\psi_\lambda$ and wavelets normalized according to
\iref{norm} with $r$ given by \iref{rscale}, we have the following norm equivalences (see \cite{Co,Dau,Me}):

\begin{enumerate}
\item
For Besov spaces,
\be
\|f\|_{\dot B^t_{q,b}}\sim \(\sum_{j\in\Z} (\sum_{|\lambda|=j}  |d_\lambda|^q)^{b/ q}\)^{1/ b},
\label{besov}
\ee
with the standard modification when $q=\infty$ or $b=\infty$.
\item
For Triebel-Lizorkin spaces,
\be
\|f\|_{\dot F^t_{q,b}}\sim \left \| \(\sum_{\lambda\in\nabla} |d_\lambda \Chi_\lambda|^b\)^{1/ b}\right \|_{L^q},
\label{triebel}
\ee
where $\Chi_\lambda=2^{dj/q}\Chi(2^j\cdot-k)$ with
$\Chi=\Chi_{[0,1]^d}$ for $\lambda\sim(j,k)$. When $b=\infty$, $\(\sum_{\lambda\in\nabla} |d_\lambda \Chi_\lambda|^b\)^{1/ b}$
should be replaced by $\sup_{j\in\Z}\sum_{|\lambda|=j} |d_\lambda \Chi_\lambda|$.
\item
For $BMO$,
\be
\|f\|_{BMO}\sim \max_{\lambda\in\nabla}  \(2^{d|\lambda|}\sum_{\mu \subset \lambda} |d_\mu|^2 2^{-d|\mu|}\)^{1/2},
\label{bmo}
\ee
where by definition $\mu \subset \lambda$ means that $2^{-j(\mu)}([0,1]^d+k(\mu)) \subset 2^{-j(\lambda)}([0,1]^d+k(\lambda))$.
\end{enumerate}
Note that due to the discretization of the scale-space index $d_\lambda$,
the above equivalent norms do not exactly satisfy the scaling relation \iref{scaling}.
These norm equivalences readily imply that $(\psi_\lambda)_{\lambda\in\Lambda}$ is
an unconditional basis for such spaces. Note that there exists no simple wavelet characterization
of Lorentz spaces $L^{q,b}$ when $b\neq q$. However the unconditionality of
$(\psi_\lambda)_{\lambda\in\Lambda}$ in such spaces follows by interpolation of Lebesgue spaces
for any $1<b,q<\infty$.

We now need to discuss the validity of Assumptions 1 and 2, for such choices of spaces.
We first discuss Assumption 2 which is only concerned with the space $X$.
Since we assumed here that $X$ is of Besov type $\dot B^{s}_{p,a}$ or Triebel-Lizorkin type $\dot F^s_{p,a}$,
we may use equivalent norms given by \iref{besov} and \iref{triebel}. Therefore, the $X$ norm of
$f_n=\sum_{\lambda\in\nabla} c_{\lambda,n}\psi_\lambda$ is either equivalent to
$$
\(\sum_{j\in\Z} (\sum_{|\lambda|=j}  |c_{\lambda,n}|^p)^{a/ p}\)^{1/ a},
$$
or
$$
\left \| \(\sum_{\lambda\in\nabla} |c_{\lambda,n} \Chi_\lambda|^a\)^{1/ a}\right \|_{L^p}.
$$
In both cases, we may invoke Fatou's lemma to conclude that for the limit sequence $(c_\lambda)$, we have
$$
\(\sum_{j\in\Z} (\sum_{|\lambda|=j}  |c_{\lambda}|^p)^{a/ p}\)^{1/ a} \leq \liminf_{n\to +\infty}\(\sum_{j\in\Z} (\sum_{|\lambda|=j}  |c_{\lambda,n}|^p)^{a/ p}\)^{1/ a},
$$
and
$$
\left \| \(\sum_{\lambda\in\nabla} |c_{\lambda} \Chi_\lambda|^a\)^{1/ a}\right \|_{L^p}
\leq   \liminf_{n\to +\infty}\left \| \(\sum_{\lambda\in\nabla} |c_{\lambda,n} \Chi_\lambda|^a\)^{1/ a}\right \|_{L^p}.
$$
Therefore, Assumption 2 holds for all Besov and Triebel-Lizorkin spaces.
\nl

We next discuss Assumption 1, for some specific examples of pairs $X$ and $Y$ which satisfy the
critical embedding property. The study of the nonlinear projector $Q_M$ is an
important chapter of approximation theory. The process of approximating
a function
$$
f=\sum_{\lambda\in\nabla} d_\lambda\psi_\lambda,
$$
by a function of the form
$$
\sum_{\lambda\in E_M} c_\lambda\psi_\lambda,
$$
with $\#(E_M)\leq M$ is sometimes called {\it best $M$-term approximation},
and has been studied extensively. The
most natural choice is to take for $E_M$ the indices corresponding to
the largest coefficients $|d_\lambda|$ and to set $c_\lambda=d_\lambda$, which corresponds to
our definition of $Q_M$. However as already mentioned in Remark \ref{remarknonlin},
other more relevant choices could be used if necessary for proving the validity
of Assumption 1 for certain pairs $(X,Y)$ and a specific instance will be mentioned below.

The study of the convergence of $Q_Mf$ towards $f$ is particularly elementary
in the case where $X=\dot B^s_{p,p}$ and $Y=\dot B^t_{q,q}$,
with $\frac 1 p-\frac 1 q=\frac {s-t} d$. Indeed, according to \iref{besov}, we have for such spaces
$$
\|f\|_{\dot B^s_{p,p}}\sim \|(d_\lambda)_{\lambda\in\nabla}\|_{\ell^p}\;\;{\rm and}\;\; \|f\|_{\dot B^t_{q,q}}\sim \|(d_\lambda)_{\lambda\in\nabla}\|_{\ell^q},
$$
and therefore for any $f\in X$, using the
decreasing rearrangement $(d_m)_{m>0}$ of the $|d_\lambda|$, we obtain
$$
\begin{array}{ll}
\|f-Q_Mf\|_{\dot B^t_{q,q}}&Ê\sim \(\sum_{\lambda\notin E_M} |d_\lambda|^q\)^{1/q} \\
& = \(\sum_{m>M} |d_m|^{Êq}\)^{1/q}  \\
& \leq |d_M|^{1-p/q}\(\sum_{m>M} |d_m|^p\)^{1/q} \\
& \leq \(M^{-1}\sum_{m=1}^M |d_m|^p\)^{1/p-1/q} \(\sum_{m>M} |d_m|^p\)^{1/q} \\
& \leq M^{-(1/p-1/q)} \(\sum_{m>0} |d_m|^p\)^{1/p}  \\
& \leq M^{-\frac {s-t} d}  \|(d_\lambda)_{\lambda\in\nabla}\|_{\ell^p} \sim M^{-\frac {s-t} d}\|f\|_{\dot B^s_{p,p}}.
\end{array}
$$
We have thus proved that
\be
\sup_{\|f\|_{\dot B^s_{p,p}} \leq 1}\|f-Q_Mf\|_{\dot B^t_{q,q}}\leq CM^{-\sigma},\;\; \sigma:=\frac {s-t} d>0,
\label{nonlinbesov}
\ee
which shows that Assumption 1 holds in such a case.

For other choices of $X$ and $Y$, the study of best $M$-term approximation is more
involved and we just describe the available results without proof.

The case
of the embedding of the Besov space $X=\dot B^{s}_{p,p}$ into
the Lebesgue space $Y=L^q$, with $q<\infty$ and $\frac 1 p-\frac 1 q=\frac {s} d>0$
has first been treated in \cite {DJP} - see also \cite{Co} and \cite{De} - where it was
proved that
\be
\sup_{\|f\|_{\dot B^s_{p,p}}  \leq 1}\|f-Q_Mf\|_{L^q}\leq CM^{-\sigma},\;\; \sigma:=\frac {s} d>0.
\label{nonlinbesleb}
\ee
Therefore Assumption 1 also holds in such a case.
Note that when $q\leq 2$, one has continuous embedding of $\dot B^0_{q,q}$ in $L^q$
and therefore \iref{nonlinbesleb} may be viewed as a consequence of \iref{nonlinbesov},
however this is no more the case when $2\leq q<\infty$, yet \iref{nonlinbesleb} still holds.

A finer result, that may be obtained by interpolation techniques, states
that, with the same relations between $p$ and $q$, the
Besov space $\dot B^{s}_{p,q}$ - which is strictly larger than $\dot B^s_{p,p}$ is continuously embedded in $L^q$,
and one may therefore ask if Assumption 1 is still valid in such a case. A positive
answer was given in \cite{K} for the more general embedding of
$X=\dot B^{s}_{p,q}$ into $Y=\dot F^t_{q,b}$ with $\frac 1 p-\frac 1 q=\frac {s-t} d$,
where $b\in ]0,\infty]$ is arbitrary: we have the convergence estimate
\be
\sup_{\|f\|_{\dot B^s_{p,q}} \leq 1}\|f-Q_Mf\|_{\dot F^t_{q,b}}\leq CM^{-\sigma},\;\; \sigma:=\frac {s-t} d>0,
\label{nonlinbesovtrib}
\ee
and therefore Assumption 1
is again valid. Note that $L^q=\dot F^0_{q,2}$ is a particular case.

For Besov spaces, the critical embedding of
$X=\dot B^s_{p,a}$ into $Y=\dot B^t_{q,b}$
with $\frac 1 p-\frac 1 q=\frac {s-t} d$  is known
to hold whenever $a\leq b$ (it is an immediate consequence
of the norm equivalence \iref{besov}). The study of
best $M$-term approximation in this context was
done in \cite{K}, where the following result was proved:
there exists a nonlinear projector $Q_M$ of the form \iref{nonlinear}, such that
when $\frac 1 a-\frac 1 b \geq \frac {s-t} d$, one has
\be
\sup_{\|f\|_{\dot B^s_{p,a}} \leq 1}\|f-Q_Mf\|_{\dot B^t_{q,b}}\leq CM^{-\sigma},\;\; \sigma:=\frac {s-t} d>0.
\label{nonlinbesov2}
\ee
The set $E_M(f)$ used in the definition of $Q_M$ is however not generally based
on picking the $M$ largest $|d_\lambda|$, which is not a problem
for our purposes as already mentioned in Remark \ref{remarknonlin}.
Therefore, Assumption 1 is valid for such pairs.

In this last example, the restriction $\frac 1 a-\frac 1 b \geq \frac {s-t} d$
is stronger than $a\leq b$ which is sufficient for the critical embedding.
However, we may still obtain the validity of Assumption 1
when $a<b$ by a general trick which we shall
re-use further: introduce an auxiliary space $Z$ with continuous embedding
\be
X\subset Z \subset Y,
\label{triple}
\ee
such that Assumption 1 either holds for the embedding between
$X$ and $Z$, or between $Z$ and $Y$, which immediately
implies the validity of Assumption 1 between $X$ and $Y$.
In the present case we take
$$
Z=\dot B^{\t s}_{\t p, a}\;\;{\rm with} \;\; t<\t s <s,\;\; \frac 1 p-\frac 1 {\t p}=\frac {s-\t s} d.
$$
The continuous embeddings \iref{triple} clearly hold.
In addition, when $\t s$ sufficiently close to $t$, we have that
$\frac 1 a-\frac 1 b \geq \frac {\t s-t} d$, so that
Assumption 1 is valid for the pair $(Z,Y)$
according to \iref{nonlinbesov2}, and thus also for $(X,Y)$.

\begin{remark}
\label{rembes}
It is not difficult to check that Assumption 1 does not hold
for the embedding of $X=\dot B^s_{p,a}$ into $Y=\dot B^t_{q,a}$,
and we also conjecture that the profile decomposition does
not generally exist for such an embedding. As an example,
consider $a=\infty$, and a sequence $(u_n)_{n>0}$ obtained
by piling up one wavelet at each scale $j=0,\cdots,n$ at position $k=0$:
$$
u_n=\sum_{j=0}^n 2^{rj}\psi(2^j\cdot).
$$
All wavelets in $u_n$ contribute equally to the $X$ and $Y$ norm
(which is equivalent to the supremum of the coefficients, equal to $1$)
and the extraction of profiles with asymptotically orthogonal scale-space
localization seems impossible.
\end{remark}

The above trick based on the intermediate space $Z$
may be used to prove Assumption 1 for other types of critical embeddings:
\begin{itemize}
\item
Embedding of Besov spaces into $BMO$:
$$
\dot B^s_{p,p} \subset BMO,\;\; s=\frac d p>0,
$$
which includes as a particular case the well known embedding
$\dot H^{d/2}\subset BMO$, and may be easily proved
from the wavelet characterization \iref{besov} and \iref{bmo}.
Choosing $Z=\dot B^{\t s}_{\t p, \t p}$ for any $0<\t s <s$
and $\t p$ such that $\t s =\frac d {\t p}$,
we clearly have the continuous embeddings \iref{triple}.
In addition, Assumption 1 is valid for the pair $(X,Z)$
according to \iref{nonlinbesov}, and thus also for $(X,Y)$.
\item
Embedding of Besov spaces into Lorentz spaces:
$$
\dot B^s_{p,a} \subset L^{q,b},\;\; \frac 1 p-\frac 1 q=\frac {s} d>0,
$$
which is valid for any $a\leq b$. If $a<b$, we may introduce for any $0<\t s < s$
$$
Z= \dot B^{\t s}_{\t p,b},\;\; \frac 1 {\t p}-\frac 1 q=\frac {\t s} d>0,
$$
so that we have the continuous embeddings \iref{triple}. In addition, we
have already proved that Assumption 1 holds
for the pair $(X,Z)$. It therefore holds for the pair $(X,Y)$.
One may easily check
that Assumption 1 does not hold
for the embedding of $X=\dot B^s_{p,a}$ into $Y=L^{q,a}$,
and conjecture that  the profile decomposition does
not generally exist for such an embedding, by an argument
analogous to the one in Remark \ref{rembes}.
\item
Embedding of Triebel-Lizorkin spaces into Triebel-Lizorkin or Besov spaces:
for any $a,b>0$, consider $X=\dot F^s_{p,a}$
and $Y=\dot F^t_{q,b}$ with $\frac 1 p-\frac 1 q=\frac {s-t} d$.
It is known - see \cite{triebel} - that $X$ is continuously embedded into
$$
Z=\dot B^{\t s}_{\t p,\t p},\;\; \t s <s\;\; {\rm and}\;\; \frac 1 p-\frac 1 {\t p}=\frac {s-\t s} d.
$$
If we assume $t<\t s< s$, we have the continuous embeddings \iref{triple}.
Moreover, we have already proved that Assumption 1 holds for the pair $(Z,Y)$.
It therefore holds for the pair $(X,Y)$.
The same type of reasoning allows to prove Assumption 1 for the embedding of $X=\dot F^s_{p,a}$
into $Y=\dot B^t_{q,b}$ with $b> p$.
\end{itemize}

\begin{remark}
It is easily seen that if $X$ and $Y$ are a pair of spaces such that both Assumptions 1 and 2 hold
for a certain wavelet basis $(\psi_\lambda)$, then
the corresponding vector fields spaces $(X)^d$ and $(Y)^d$
also satisfy the same assumptions for the
vector valued wavelet basis
$$
\psi_{\lambda,i}:=\psi_\lambda e_i, \;\; \lambda\in\nabla, \;i=1,\cdots,d,
$$
where $e_i=(0,\cdots,0,1,0,\cdots,0)$ is the canonical basis vector, and
the wavelet coefficients $d_\lambda$ are defined accordingly as vectors.
\end{remark}

\section{Stability of the decomposition}

We finally want to show that the decomposition is stable
in the sense that the sum of the $\|\phi^l\|_X$ raised
to an appropriate power remains bounded.
In our discussion, we distinguish between the cases
where $X$ is a Besov or Triebel-Lizorkin space.
We first address the Besov case.

\begin{theorem}
Assume that $X=\dot B^s_{p,a}$ with $1\leq p<\infty$ and $1\leq a\leq \infty$.
We then have
\be
\| (\|\phi^l\|_X)_{l>0} \|_{\ell^\tau}\leq CK, \;\;  \tau:=\max\{p,a\}.
\label{stabbesov}
\ee
where $C$ is a constant that only depends on $X$ and on the
choice of the wavelet basis and where $K:=\sup_{n\geq 0} \|u_n\|_X$.
\end{theorem}

\noindent
{\bf Proof:} Fix an arbitrary $L>0$ and let $M$ be such
that $L\leq L(M)$ as in Step 4 of \S 2. For $l=1,\cdots,L$,
we recall the approximate
profiles
$$
\phi^{l,M}_{\lambda_l(n)}=\sum_{m\in E(l,M)}d_m\psi_{\lambda(m,n)},
$$
and we also define
$$
\phi^{l,M,n}:=\sum_{m\in E(l,M)}d_{m,n}\psi_{\lambda(m,n)},
$$
which are disjoint part of the wavelet expansion of $u_n$. More
precisely, we have
$$
u_n=\sum_{l=1}^L \phi^{l,M,n}+\sum_{m>M} d_{m,n}\psi_{\lambda(m,n)}.
$$
We next claim that if $E_1,\cdots,E_L$ are disjoints finite sets
in $\nabla$, then for any coefficient sequence $(d_\lambda)$, one has
\be
\( \sum_{l=1}^{L} \| \sum_{\lambda\in E_l} d_\lambda\psi_\lambda\|_X^\tau\)^{1/\tau}  \leq C\left \|\sum_{l=1}^{L}  \sum_{\lambda\in E_l} d_\lambda\psi_\lambda\right \|_X,
\label{claimbesov}
\ee
where $C$ is a constant that only depends on $X$ and on the choice
of the wavelet basis, and with the standard modification of the
sum to the power $1/\tau$ by a supremum in the left hand side when $\tau=\infty$.

Before proving this claim, we first show that it leads to the conclusion of the proof. Indeed, for $l=1,\cdots,L$,
the functions $\phi^{l,M,n}$ are linear combinations of wavelets with indices in disjoint
finite sets $E_1,\cdots,E_L$
(that vary with $n$), and therefore according to \iref{claimbesov}, when $\tau<\infty$,
$$
\(\sum_{l=1}^L  \|\phi^{l,M,n}\|_X^\tau\)^{1/\tau} \leq C
\left \| \sum_{l=1}^L \phi^{l,M,n}\right \|_X.
$$
Using the unconditionality inequality \iref{uncond}, we thus find that for all $n>0$
$$
\(\sum_{l=1}^L  \|\phi^{l,M,n}\|_X^\tau\)^{1/\tau}\leq CK,
$$
up to a multiplication of the constant $C$ by $D$.
Since $\|\phi^{l,M}_{\lambda_l(n)}-\phi^{l,M,n}\|_X\to 0$ as $n\to \infty$, it follows that for any
$\e>0$ we have
$$
\(\sum_{l=1}^L  \|\phi^{l,M}_{\lambda_l(n)}\|_X^\tau\)^{1/\tau}\leq CK+\e,
$$
for $n$ sufficiently large. By the scaling invariance \iref{scaling} we thus find that
$$
\(\sum_{l=1}^L  \|\phi^{l,M}\|_X^\tau\)^{1/\tau} \leq CK.
$$
Letting $M$ go to $+\infty$, we obtain the same inequality for
the exact profiles
$$
\(\sum_{l=1}^L  \|\phi^{l}\|_X^\tau\)^{1/\tau} \leq CK,
$$
and we thus conclude that \iref{stabbesov} holds, by letting $L\to +\infty$. The case $\tau=\infty$
is treated in an exact similar way, replacing the sum to the power $1/\tau$ by a supremum.

It remains to prove \iref{claimbesov}. We actually claim that
this property holds with constant $C=1$ if we take for $\|\cdot\|_X$ the
equivalent norm given by \iref{besov}. This is obvious when $p=a=\tau$
since this equivalent norm is then simply the $\ell^\tau$ norm of
the wavelet coefficients. When $p\neq a$, we distinguish between
the cases $p<a$ and $p>a$. We denote by
$$
E_{j,l}:=\{\lambda\in E_l\; ; \; |\lambda|=j\},
$$
so that
$$
E_l=\cup_{j\in\Z} E_{l,j}.
$$
First consider the case $\tau=a>p$. We then have, when $a<\infty$,
$$
\begin{array}{ll}
\left \|\sum_{l=1}^{L}  \sum_{\lambda\in E_l} d_\lambda\psi_\lambda\right \|_X^a
& = \sum_{j\in\ZZ}\(\sum_{l=1}^L\sum_{\lambda\in E_{l,j}}|d_\lambda|^p\)^{a/p}\\
&    \geq \sum_{j\in\ZZ}\sum_{l=1}^L\(\sum_{\lambda\in E_{l,j}}|d_\lambda|^p\)^{a/p}\\
& = \sum_{l=1}^L \sum_{j\in\ZZ}\(\sum_{\lambda\in E_{l,j}}|d_\lambda|^p\)^{a/p}\\
&= \sum_{l=1}^L \| \sum_{\lambda\in E_l} d_\lambda\psi_\lambda \|_X^a,
\end{array}
$$
where for the inequality we have simply used the fact that $a/p>1$.
Therefore \iref{claimbesov} holds.  When $a=\infty$, we obtain the same result
by writing
$$
\begin{array}{ll}
\left \|\sum_{l=1}^{L}  \sum_{\lambda\in E_l} d_\lambda\psi_\lambda\right \|_X
& = \sup_{j\in\ZZ}\(\sum_{l=1}^L\sum_{\lambda\in E_{l,j}}|d_\lambda|^p\)^{1/p}\\
&    \geq \sup_{j\in\ZZ}\(\sup_{l\leq L}\sum_{\lambda\in E_{l,j}}|d_\lambda|^p\)^{1/p}\\
& = \sup_{l\leq L}  \sup_{j\in\ZZ}\(\sum_{\lambda\in E_{l,j}}|d_\lambda|^p\)^{1/p}\\
&= \sup_{l\leq L} \| \sum_{\lambda\in E_l} d_\lambda\psi_\lambda \|_X.
\end{array}
$$
We next consider the case $\tau=p>a$.  To treat this case where $p < \infty$ by hypothesis, we introduce the notation
$$
b_{j,l}:=\(\sum_{\lambda\in E_{l,j}}|d_\lambda|^p\)^{a/p},
$$
and we remark that \iref{claimbesov} is then
equivalent to
$$
\(\sum_l \( \sum_j b_{j,l}\)^{p/a}\)^{a/p} \leq  \sum_j \( \sum_l |b_{j,l}|^{p/a}\)^{a/p},
$$
which trivially holds by applying the triangle inequality in $\ell^{p/a}$.  \hfill $\Box$
\nl

Our last result addresses the Triebel-Lizorkin case.

\begin{theorem}
Assume that $X=\dot F^s_{p,a}$ with $1\leq p<\infty$ and $1\leq a\leq \infty$.
We then have
\be
\| (\|\phi^l\|_X)_{l>0} \|_{\ell^p}\leq CK.
\label{stabTriebel}
\ee
where $C$ is a constant that only depends on $X$ and on the
choice of the wavelet basis and where $K:=\sup_{n\geq 0} \|u_n\|_X$.
\end{theorem}

\noindent
{\bf Proof:} We only give the proof in the case $a<\infty$,
the case $a=\infty$ being treated by the same type
of arguments up to notational changes.
Fix an arbitrary $L>0$ and let $M$ be such
that $L\leq L(M)$ as in Step 4 of \S 2. By the unconditionality
\iref{uncond} of the wavelet basis with respect to $X$, we first observe that
$$
\left\|\sum_{m=1}^M d_{m,n}\psi_{\lambda(m,n)}\right \|_X\leq DK.
$$
It follows that for any $\e>0$, we have
$$
\left\|\sum_{m=1}^M d_{m}\psi_{\lambda(m,n)}\right \|_X\leq DK+\e.
$$
for $n$ sufficiently large. Recall that the sum inside the norm
may be rewritten in terms of the approximate profiles:
$$
\sum_{m=1}^M d_{m}\psi_{\lambda(m,n)}=\sum_{l=1}^L\sum_{m\in E(l,M)} d_{m}\psi_{\lambda(m,n)}=\sum_{l=1}^L \phi^{l,M}_{\lambda_l(n)}.
$$
We associate to the approximate profile $\phi^{l,M}$ a piecewise constant function $\Chi^{l,M}$ defined by
$$
\Chi^{l,M}_{\lambda_l(n)}:=\sum_{m\in E(l,M)} |d_{m}\Chi_{\lambda(m,n)}|^a,
$$
where $\Chi_{\lambda}=2^{dj/p}\Chi(2^j\cdot - k)$ with $\Chi=\Chi_{[0,1]^d}$ for $\lambda\sim (j,k)$. Thus
according to the wavelet characterization \iref{triebel} of Triebel-Lizorkin spaces, we have
$$
c \int_{\RR^d}|\Chi^{l,M}_{\lambda_l(n)}(x)|^{p/a}dx\leq \|\phi^{l,M}_{\lambda_l(n)}\|_X^p\leq C\int_{\RR^d}|\Chi^{l,M}_{\lambda_l(n)}(x)|^{p/a}dx,
$$
as well as
$$
c\int_{\RR^d} |\sum_{l=1}^L\Chi^{l,M}_{\lambda_l(n)} (x)|^{p/a} dx\leq  \left\|\sum_{m=1}^M d_{m,n}\psi_{\lambda(m,n)}\right \|_X^p
\leq C\int_{\RR^d} |\sum_{l=1}^L\Chi^{l,M}_{\lambda_l(n)}(x) |^{p/a}dx ,
$$
where $0<c\leq C$ only depend on the choice of the wavelet basis. In the case where $a\leq p$, we obviously have
$$
\sum_{l=1}^L\int_{\RR^d}|\Chi^{l,M}_{\lambda_l(n)}(x)|^{p/a}dx\leq \int_{\RR^d} |\sum_{l=1}^L\Chi^{l,M}_{\lambda_l(n)} (x)|^{p/a}dx.
$$
It therefore follows that
$$
\sum_{l=1}^L \|\phi^{l,M}\|_X^p=\sum_{l=1}^L \|\phi^{l,M}_{\lambda_l(n)}\|_X^p\leq \frac C c \left\|\sum_{m=1}^M d_{m,n}\psi_{\lambda(m,n)}\right \|_X^p
\leq  \frac C c (DK+\e)^p.
$$
Since this holds for any $\e>0$, and $L>0$ and $M$ such
that $L\leq L(M)$, we therefore obtain \iref{stabTriebel}
by a limiting argument, up to renaming $(C/c)^{1/p}D$ into $C$.

In order to reach the same conclusion in the case $p<a$, we need to exploit the
``asymptotic orthogonality'' of the scales $\lambda_l(n)$ as expressed by \iref{orthoscale} in the statement
of Theorem \iref{maintheo}. For this purpose, let us define
$$
\Omega_n:={\rm Supp}\(\sum_{l=1}^L\Chi^{l,M}_{\lambda_l(n)}\)=\cup_{l=1}^L {\rm Supp}(\Chi^{l,M}_{\lambda_l(n)}).
$$
For any $x\in \Omega_n$, we denote by $l^*$ the number in $\{1,\cdots,L\}$ such that
$$
\Chi^{l^*,M}_{\lambda_{l^*}(n)}(x)=\max_{l=1,\cdots,L} \Chi^{l,M}_{\lambda_l(n)}(x).
$$
Note that $l^*$ depends both of $x$ and $n$.
We claim that a consequence of \iref{orthoscale} is that the function $\Chi^{l^*,M}_{\lambda_{l^*}(n)}$
tends to dominate all other $\Chi^{l,M}_{\lambda_l(n)}$ at the point $x$ as $n\to +\infty$ in the following uniform sense:
\be
\lim_{n\to +\infty} \min_{x\in\Omega_n} \frac{\Chi^{l^*,M}_{\lambda_{l^*}(n)}(x)}{\sum_{l=1}^L \Chi^{l,M}_{\lambda_l(n)}(x)
-\Chi^{l^*,M}_{\lambda_{l^*}(n)}(x)}  =+\infty.
\label{claimtriebel}
\ee
Before proving this claim, let us show how it leads us to the conclusion of the theorem. We observe that
\iref{claimtriebel} also means that
$|\Chi^{l^*,M}_{\lambda_{l^*}(n)}|^{p/a}$
tends to dominate all other $|\Chi^{l,M}_{\lambda_l(n)}|^{p/a}$ at the point $x$ as $n\to +\infty$.
Therefore, for any $\e>0$, we
have for $n$ large enough
$$
\sum_{l=1}^L| \Chi^{l,M}_{\lambda_l(n)}(x)|^{p/a}
\leq (1+\e) |\Chi^{l^*,M}_{\lambda_{l^*}(n)}(x)|^{p/a} \leq (1+\e) |\sum_{l=1}^L \Chi^{l,M}_{\lambda_l(n)}(x)|^{p/a},
$$
for all $x\in\Omega_n$, and thus
$$
\sum_{l=1}^L\int_{\RR^d}|\Chi^{l,M}_{\lambda_l(n)}(x)|^{p/a}dx\leq (1+\e) \int_{\RR^d} |\sum_{l=1}^L\Chi^{l,M}_{\lambda_l(n)} (x)|^{p/a}dx.
$$
We may then conclude the proof as in the case $a\leq p$.

It remains to prove \iref{claimtriebel}. Our
first observation is that
the asymptotic orthogonality of the scales $\lambda_l(n)$ expressed by \iref{orthoscale},
shows that for a given $x$ the profile scales $|\lambda_l(n)|$ for those $l\in \{1,\cdots,L\}$
such that $x\in {\rm Supp}( \Chi^{l,M}_{\lambda_l(n)})$ tend to get far apart as $n$ grows.
Indeed, these $\lambda_l(n)$ do not get far apart in space since the supports
of $\Chi^{l,M}_{\lambda_l(n)}$ all contain the same point $x$.

We introduce $l_*$ the number that
maximizes $|\lambda_l(n)|$ among all those $l\in \{1,\cdots,L\}$ such that
$x\in {\rm Supp}( \Chi^{l,M}_{\lambda_l(n)})$. Similar to $l^*$,
the number $l_*$ depends on both $x$ and $n$. From the previous
observation, we know that for
any arbitrarily large $B>0$, there exists $n_0$ such that
for all $n\geq n_0$, we have
\be
|\lambda_{l_*}(n)|\geq |\lambda_l(n)|+B,
\label{lB}
\ee
for all $l\in \{1,\cdots,L\}$ such that $x\in {\rm Supp}( \Chi^{l,M}_{\lambda_l(n)})$
and $l\neq l_*$. Moreover we may choose this $n_0$ independent of the selected point
$x$ for the same $B>0$.

We claim that as $n$ grows $\Chi^{l_*,M}_{\lambda_{l_*}(n)}$
tends to dominate all other $\Chi^{l,M}_{\lambda_l(n)}$ at the point $x$ as $n\to +\infty$,
in the sense that
\be
\lim_{n\to +\infty} \min_{x\in\Omega_n} \frac{\Chi^{l_*,M}_{\lambda_{l_*}(n)}(x)}{\sum_{l=1}^L \Chi^{l,M}_{\lambda_l(n)}(x)
-\Chi^{l_*,M}_{\lambda_{l_*}(n)}(x)}  =+\infty.
\label{claimtriebel2}
\ee
This clearly implies \iref{claimtriebel} (and shows that $l^*=l_*$ for $n$ large enough).

In order to prove \iref{claimtriebel2}, we observe that if $x\in {\rm Supp}( \Chi^{l,M}_{\lambda_l(n)})$
for some $l\in \{1,\cdots,L\}$, we may then frame $\Chi^{l,M}_{\lambda_l(n)}(x)$ according to
$$
c Ê2^{\frac{adj_l(n)}p }\leq  \Chi^{l,M}_{\lambda_l(n)}(x)\leq CÊ2^{\frac{adJ_l(n)}p }
$$
where
$$
j_l(n):=\min_{m\in E(l,M)} |\lambda(m,n)|\;\;{\rm and}\;\; J_l(n):=\max_{m\in E(l,M)} |\lambda(m,n)|,
$$
and where
$$
c:=|d_M|^a\;\;{\rm and} \;\; C:=\sum_{m=1}^M |d_m|^a.
$$
The constants $c$ and $C$ of course depend on $L$ and $M$
which are fixed at that stage.
Note that from the construction of the profile there exists $A>0$ (that also depends on $L$ and $M$)
such that for all $l\in \{1,\cdots,L\}$
$$
|\lambda_l(n)|-A \leq j_l(n)\leq J_l(n)\leq |\lambda_l(n)|+A,
$$
and therefore, up to a modification in the constants $c$ and $C$ we may write
$$
c Ê2^{\frac{ad|\lambda_l(n)|}p}\leq  \Chi^{l,M}_{\lambda_l(n)}(x)\leq CÊ2^{\frac{ad|\lambda_l(n)|}p }.
$$
Combining this observation with \iref{lB}, we easily obtain \iref{claimtriebel2}. \hfill $\Box$

\vbox{\noindent Hajer Bahouri\\
Centre de Math\'ematiques - Facult\'e de Sciences et TechnologieÊ\\
Universit\'e Paris XII - Val de Marne\\
61, avenue du G\'en\'eral de Gaulle\\
94010 Creteil Cedex, France \\
e--mail: {\tt hbahouri@math.cnrs.fr}}

\bigskip

\vbox{\noindent Albert Cohen\\
Laboratoire Jacques-Louis Lions\\
Universit\'e Pierre et Marie Curie\\
175 Rue du Chevaleret, 75013 Paris\\
France\\
e--mail: {\tt cohen@ann.jussieu.fr}\\
www: {\tt http://www.ann.jussieu.fr/$\sim$cohen}\\
Tel: 33-1-44277195,  Fax: 33-1-44277200}

\bigskip

\vbox{\noindent Gabriel Koch\\
Mathematical Institute\\
24-29 St Giles'\\
Oxford, OX1 3LB,
England \\
e--mail: {\tt koch@maths.ox.ac.uk}}

\end{document}